\documentclass[12pt]{amsart}
\usepackage{latexsym}
\usepackage{amsmath}
\usepackage{amsfonts}
\usepackage{amssymb}
\usepackage{amsthm}
\theoremstyle{plain}
\newtheorem{theorem}{Theorem}[section]
\newtheorem{lemma}[theorem]{Lemma}
\newtheorem{proposition}[theorem]{Proposition}

\newtheorem*{theorem*}{Theorem}

\theoremstyle{definition}
\newtheorem{definition}[theorem]{Definition}

\theoremstyle{remark}

\newcommand{\E}{\mathbb{E}}

\newcommand{\N}{\mathbb{N}}
\newcommand{\C}{\mathbb{C}}

\newcommand{\R}{\mathbb{R}}
\newcommand{\T}{\mathbb{T}}
\newcommand{\Z}{\mathbb{Z}}

\newcommand{\e}{\varepsilon}

\newcommand{\norm}[1]{\lVert #1\rVert}

\begin{document}
\title{Uniformity  in the polynomial Wiener-Wintner theorem}
\author{Nikos Frantzikinakis}
\begin{abstract}
In 1993, E. Lesigne proved a polynomial extension of the
Wiener-Wintner ergodic theorem
 and asked two questions: does this result have a uniform counterpart
 and can an assumption of
 total ergodicity  be replaced by ergodicity?
 The purpose of this article is to answer these questions,
 the first one positively and the second one negatively.
\end{abstract}

\address{1 Einstein Drive,
Institute for Advanced Study, Princeton, NJ 08540}
\maketitle

\section{Introduction}
The Wiener-Wintner ergodic theorem~\cite{W} asserts that if
$(X,\mathcal{B},\mu,T)$ is an ergodic measure preserving system
and
$f\in L^1(\mu)$, then  for a.e. $x\in X$  the limit
\begin{equation}\label{E:WW}
\lim_{N\to\infty}\frac{1}{N}\sum_{n=1}^N e(n\alpha)\cdot f(T^nx)
\end{equation}
exists for every $\alpha\in\R$, where $e(t)=e^{2\pi i t}$. What
makes this a nontrivial strengthening of the Birkhoff ergodic
theorem is that the set of full measure for which we have
convergence does not depend on the choice of $\alpha\in \R$.

 Two different proofs of this theorem are based on
 the following results which are of interest on their own:

$(i)$ For ergodic systems,  if $f\in L^1(\mu)$ and
$f\in\mathcal{E}_1(T)^\bot$, then for a.e. $x\in X$ we have that
for every $\alpha\in\R$
\begin{equation}\label{E:WW1}
\lim_{N\to\infty} \frac{1}{N} \sum_{n=1}^N e(n\alpha)\cdot
f(T^nx)=0,
\end{equation}
where $\mathcal{E}_1(T)$ is the set of eigenfunctions of $T$.

 $(ii)$ For ergodic systems, if $f\in L^1(\mu)$ and $f\in\mathcal{E}_1(T)^\bot$, then for a.e.
$x\in X$ we have
\begin{equation}\label{E:WW2}
\lim_{N\to\infty} \sup_{\alpha\in\R}\Big|\frac{1}{N} \sum_{n=1}^N
e(n\alpha)\cdot f(T^nx)\Big|=0.
\end{equation}
Condition $(ii)$  is  of course stronger than $(i)$.

 A natural
problem is to find conditions that will allow us to replace the
linear polynomials $P(n)=n\alpha$,  in \eqref{E:WW},
\eqref{E:WW1}, and \eqref{E:WW2}, with higher degree polynomials.
This study was initiated by Lesigne in \cite{L0} and \cite{L}. To
state his main result we need the notion of level $k$
quasi-eigenfunction, introduced by Abramov\footnote{Abramov
studied in~\cite{A} totally ergodic systems  for which the
subspace spanned by elements of $\bigcup_{k\in\N}
\mathcal{E}_k(T)$ is dense in $L^2(\mu)$. He showed that any such
system is isomorphic to a unipotent affine transformation on a
compact Abelian group.}:
\begin{definition}
If $(X,\mathcal{X},\mu,T)$ is an ergodic system we let
$\mathcal{E}_0(T)$ denote the set of eigenvalues of $T$ and for
$k\in\N$ we define inductively
$$
\mathcal{E}_k(T)=\{f\in L^\infty(\mu)\colon |f|=1 \text{ and
 } Tf \cdot \bar{f}\in\mathcal{E}_{k-1}(T)\}.$$

\end{definition}

The result of Lesigne is:
\begin{theorem}[{\bf Lesigne~\cite{L0}, \cite{L}}]\label{T:Lesigne}
Let $(X,\mathcal{X},\mu,T)$ be an  ergodic  system and $f\in
L^1(\mu)$.

(i) For a.e. $x\in X$, we have for every continuous $\phi\colon
\T\to \C$, and  polynomial $P$ with real coefficient that the
averages
$$
\frac{1}{N}\sum_{n=1}^N \phi(P(n))\cdot f(T^nx)
$$
converge as $N\to\infty$.

(ii) Let $k\in\N$. If $f\in\mathcal{E}_k(T)^\bot$ and the system
is totally ergodic $($meaning $T^r$ is ergodic for $r\in\N$$)$,
then for a.e. $x\in X$, we have for  every continuous $\phi\colon
\T\to \C$ , and $P\in\R_k[t]$ that
\begin{equation}\label{E:pWW}
\lim_{N\to\infty} \frac{1}{N} \sum_{n=1}^N \phi(P(n))\cdot
f(T^nx)=0,
\end{equation}
where $\R_k[t]$ denotes the set of polynomials with real
coefficients and degree at most $k$.
\end{theorem}
Motivated by the linear case where  $f\in\mathcal{E}_1(T)^\bot$
implies uniform convergence to zero  (see equation \eqref{E:WW2}),
Lesigne asked in \cite{L}   whether (ii) has a  uniform
counterpart. Furthermore, he asked whether the total ergodicity
assumption can be removed from  (ii). He remarked that this would
be the case if one could establish  that
$\mathcal{E}_k(T)^\bot=\mathcal{E}_k(T^m)^\bot$ for every $k,m\in
\N$, a result that is easily shown to be true when $k=1$. The
purpose of this paper is to answer these questions, the first
affirmatively and the second negatively. We also give an example
where one has $\mathcal{E}_k(T)^\bot\neq \mathcal{E}_k(T^2)^\bot$,
for $k\geq 2$.

The following is our main result and  gives a positive answer to
the first question of Lesigne:
\begin{theorem}\label{T:main}
Let $(X,\mathcal{X},\mu,T)$ be a totally ergodic  system, $f\in
L^1(\mu)$, and $k\in\N$. Then the following three conditions are
equivalent:

(i) $f\in\mathcal{E}_k(T)^\bot$.

(ii) For a.e. $x\in X$, we have  for every  continuous $\phi\colon
\T\to \C$, and $P\in\R_k[t]$ that
$$
\lim_{N\to\infty} \frac{1}{N} \sum_{n=1}^N \phi(P(n))\cdot
f(T^nx)=0.
$$

 (iii) For a.e. $x\in X$, we have  for every continuous
$\phi\colon \T\to \C$ that
\begin{equation}\label{E:main}
\lim_{N\to\infty} \sup_{P\in \R_k[t]}\Big|\frac{1}{N} \sum_{n=1}^N
\phi(P(n))\cdot f(T^nx)\Big|=0.
\end{equation}
\end{theorem}
We stress that in $(ii)$ the set of full measure does not depend
on the choice of the  function $\phi$ or the polynomial $P\in
\R_k[t]$, and in (iii) on the choice of the  function $\phi$.

The following result gives  a negative answer to the second
question of Lesigne:
\begin{theorem}\label{T:counterexamples}
(i) For $k\geq 2$, there exists an ergodic   system
$(X,\mathcal{B},\mu,T)$ such that $\mathcal{E}_k(T)^\bot\neq
\mathcal{E}_k(T^2)^\bot$.

(ii) For $k\geq 2$, Theorem~\ref{T:main} and  part (ii) of
Theorem~\ref{T:Lesigne} fail if we replace the total ergodicity
assumption with ergodicity.
\end{theorem}

In Section $2$ we prove Theorem~\ref{T:main} by combining the
method used in \cite{L} to prove the nonuniform result, with a new
elementary idea that enables us to get uniformity (this is
Proposition~\ref{P:key} but the main idea is contained in
Lemma~\ref{L:key}).
Finally, in Section $3$ we prove Theorem~\ref{T:counterexamples}.
We construct the system explicitly as a skew product extension of
a group rotation.

{\bf Acknowledgments}: I would like to thank J. Campbell for
bringing to my attention the questions studied in this article. I
would also like to thank B. Kra and E. Lesigne for helpful
comments.

\section{Proof of the uniform Polynomial Wiener-Wintner}

We will prove Theorem~\ref{T:main}. We start with a reduction:
\begin{proposition}\label{P:reduction}
It suffices to check the implication $(i)\Rightarrow (iii)$ in
Theorem~\ref{T:main} when $\phi(t)=e(t)$ and $f\in L^\infty(\mu)$.
\end{proposition}
\begin{proof}

Since every continuous function on $\T$ can be uniformly
approximated by trigonometric polynomials we see that it suffices
to verify \eqref{E:main} for $\phi(t)=e(mt)$, $m\in \Z$.
Furthermore,  since for every $m\in\Z$ we have $mP \in\R_{k}[t]$
whenever $P \in\R_{k}[t]$, we can restrict ourselves to the case
where $m=1$.

Suppose that for every $f\in L^\infty(\mu)$ the implication
$(i)\Rightarrow (iii)$ in Theorem~\ref{T:main} holds, we claim
that it also holds for every $f\in L^1(\mu)$. Let $k\in\N$ and
$f\in L^1(\mu)$ be such that $f\in\mathcal{E}_k(T)^\bot$. We first
notice that there exists a sequence of functions
$\{f_l\}_{l\in\N}$ such that $f_l\in L^\infty(\mu)$, $f_l\in
\mathcal{E}_k(T)^\bot$, and $\norm{f-f_l}_{L^1(\mu)}\leq 1/l$. To
see this, we choose $g_l\in L^\infty(\mu)$ such that
$\norm{f-g_l}_{L^1(\mu)}\leq 1/(2l)$. Let $\mathcal{D}$ be the
$\sigma$-algebra generated by all $h$-measurable sets  where $h\in
\mathcal{E}_k(T)$. Since $\E(f|\mathcal{D})=0$ and
$\norm{f-g_l}_{L^1(\mu)}\leq 1/(2l)$, we have
$\norm{\E(g_l|\mathcal{D})}_{L^1(\mu)}\leq 1/(2l)$. It is now easy
to check that the sequence $\{f_l\}_{l\in\N}$, where
$f_l=g_l-\E(g_l|\mathcal{D})$,  satisfies the advertised
conditions.

Let $X_0$ be the set of $x\in X$ for which the ergodic theorem
holds for functions of the form $f_l-f$, for all $l\in\N$.
Obviously $m(X_0)=1$. For every $x\in X$ we have
\begin{equation}\label{EE1}
\lim_{N\to\infty} \sup_{P\in \R_k[t]}\Big|\frac{1}{N} \sum_{n=1}^N
e(P(n))\cdot f(T^nx)\Big|\leq A_{l}+B_{l},
\end{equation}
where
$$
A_{l}=\lim_{N\to\infty} \sup_{P\in \R_k[t]}\Big|\frac{1}{N}
\sum_{n=1}^N e(P(n))\cdot (f(T^nx)-  f_l(T^nx)) \Big|
$$
and
$$
B_{l}=\lim_{N\to\infty} \sup_{P\in \R_k[t]}\Big|\frac{1}{N}
\sum_{n=1}^N e(P(n))\cdot f_l(T^nx)\Big|.
$$
For $x\in X_0$ we have
$$
A_{l}\leq \lim_{N\to\infty}\frac{1}{N} \sum_{n=1}^N |f(T^nx)-
f_l(T^nx)|=\int |f(x)-f_l(x)|\ d\mu\leq 1/l.
$$
We also know that there exists a set of full measure $X_1$ such
that $B_l=0$ for $x\in X_1$. Set $X_2=X_0\cap X_1$. Letting
$l\to\infty$ in \eqref{EE1} gives that for $x\in X_2$ we have
$$
\lim_{N\to\infty} \sup_{P\in \R_k[t]}\Big|\frac{1}{N} \sum_{n=1}^N
e(P(n))\cdot f(T^nx)\Big|=0.
$$
 Since $X_2$ has full measure the claim is proved.
\end{proof}

 The following ``uniformization'' trick contains the main idea
 needed for
Proposition~\ref{P:key}, which in turn is one of the key
ingredients needed for the proof of  Theorem~\ref{T:main}. We
prove it in order to make the argument of Proposition~\ref{P:key}
more transparent.
\begin{lemma}\label{L:key}
Let $b(n)$ be a bounded sequence of complex numbers, such that for
every $\alpha\in [0,1]$ we have
\begin{equation}\label{E:111}
\lim_{M\to\infty}\limsup_{N\to\infty} \frac{1}{N}\sum_{n=1}^{N}
\Big|\frac{1}{M} \sum_{m=1}^M e(m\alpha)\cdot b(Mn+m)\Big|=0.
\end{equation}
Then
\begin{equation}\label{E:333}
\lim_{N\to\infty} \sup_{\alpha\in [0,1]}
\Big|\frac{1}{N}\sum_{n=1}^N e(n\alpha)\cdot b(n)\Big|=0.
\end{equation}
\end{lemma}
\begin{proof}
Let $\e>0$. We can assume that $b(n)$ is bounded by $1$. Then for
every $M,N\in \N$ and $\alpha,\beta \in [0,1]$, we have
\begin{align*}
 \Big|\frac{1}{N}\sum_{n=1}^N e(n\beta)\cdot b(n)\Big|&\leq
\frac{1}{[N/M]}\sum_{n=1}^{[N/M]} \Big|\frac{1}{M} \sum_{m=1}^M
e(m\beta)\cdot b(Mn+m)\Big|+2M/N\\&\leq A_{M,N}+B_{M,N},
\end{align*}
where
$$
A_{M,N}=\frac{1}{[N/M]}\sum_{n=1}^{[N/M]} \frac{1}{M} \sum_{m=1}^M
|e(m\beta)-e(m\alpha)|
$$
and
$$
B_{M,N}=\frac{1}{[N/M]}\sum_{n=1}^{[N/M]} \Big|\frac{1}{M}
\sum_{m=1}^M e(m\alpha)\cdot b(Mn+m)\Big| +2M/N.
$$
It follows from  \eqref{E:111}, that for fixed $\alpha\in [0,1]$,
we can choose $M_\alpha$ and $N_\alpha$  such that for all
$N>N_\alpha$ we have $B_{M_\alpha,N}\leq \e/2$. We can choose a
neighborhood $V_\alpha$ of $\alpha$ such that $\sup_{\beta\in
V_\alpha,1\leq m\leq M_\alpha}|e(m\beta)-e(m\alpha)|\leq \e/2$.
Then for every $N\in \N$ and $\beta \in V_\alpha$ we have
$A_{M_a,N}\leq \e/2$. Putting these two estimates together we get
that if $N>N_\alpha$ then
 \begin{equation}\label{E:444}
\sup_{\beta\in V_\alpha} \Big|\frac{1}{N}\sum_{n=1}^N
e(n\beta)\cdot b(n)\Big|\leq \e.
\end{equation}
Doing this for every $\alpha\in [0,1]$, we  produce an open cover
$\{V_\alpha\}_{\alpha\in[0,1]}$ of $[0,1]$  and positive integers
$\{N_\alpha\}_{\alpha\in[0,1]}$, such that \eqref{E:444} holds for
$N>N_\alpha$. By compactness, there exists a finite subcover, say
$V_{\alpha_1},\ldots, V_{\alpha_l}$. Then for
$N>\max\{N_{\alpha_1},\ldots, N_{\alpha_l}\}$ we have
 $$
\sup_{\alpha\in [0,1]} \Big|\frac{1}{N}\sum_{n=1}^N
e(n\alpha)\cdot b(n)\Big|\leq \e.
$$
This proves \eqref{E:333}.
\end{proof}

\begin{proposition}\label{P:key}
Let $(X,\mathcal{X},\mu,T)$ be a totally ergodic  system, and
$f\in L^\infty(\mu)$ be such that for every $\alpha\in\R$  we have
that for a.e. $x\in X$ $($the set of full measure may depend on
the choice of $\alpha$$)$
\begin{equation}\label{E:assumption}
\lim_{N\to\infty} \sup_{P\in \R_{k-1}[t]}\Big|\frac{1}{N}
\sum_{n=1}^N  e(n^k\alpha+P(n))\cdot f(T^nx)\Big|=0.
\end{equation}
Then for a.e. $x\in X$ we have
\begin{equation}\label{E:-1}
\lim_{N\to\infty} \sup_{P\in \R_{k}[t]}\Big|\frac{1}{N}
\sum_{n=1}^N  e(P(n))\cdot f(T^nx)\Big|=0.
\end{equation}
\end{proposition}
\begin{proof}
Let $f\in L^\infty(\mu)$. Without loss of generality we can assume
that $|f(x)| \leq 1$ for every $x\in X$.

 {\bf First step.} Motivated by  Lemma~\ref{L:key} we first show
  that there exists a set of full measure $X_0$, such that
 for  $x\in X_0$,  and for every $\alpha\in\R$ we have
\begin{equation}\label{E:key}
\lim_{M\to\infty}\lim_{N\to\infty} \frac{1}{N}\sum_{n=1}^{N}
\sup_{P\in \R_{k-1}[t]}\Big|\frac{1}{M} \sum_{m=1}^M
e(m^k\alpha+P(m))\cdot f(T^{Mn+m}x)\Big|=0.
\end{equation}
To see this, for $M\in\N$  and $\alpha\in\mathbb{Q}$, we apply the
ergodic theorem for the transformation $T^M$ (we use total
ergodicity here) and the function
$$
g_{M,\alpha}(x)=\sup_{P\in \R_{k-1}[t]}\Big|\frac{1}{M}
\sum_{m=1}^M e(m^k\alpha+P(m))\cdot f(T^{m}x)\Big|.
$$
Since we  are interested in only countably many  functions and
transformations, we get a set $X_0$ of full measure, such that for
every  $x\in X_0$,  $M\in \N$, and $\alpha\in\mathbb{Q}$ we have
\begin{gather}\label{E:11}
\lim_{N\to\infty} \frac{1}{N}\sum_{n=1}^{N} \sup_{P\in
\R_{k-1}[t]}\Big|\frac{1}{M} \sum_{m=1}^M e(m^k\alpha+P(m))\cdot
f(T^{Mn+m}x)\Big|=\\ \notag \int \sup_{P\in
\R_{k-1}[t]}\Big|\frac{1}{M} \sum_{m=1}^M e(m^k\alpha+P(m))\cdot
f(T^{m}x)\Big|\ d\mu.
\end{gather}
Since the set of $\alpha\in \R$ for which \eqref{E:11} holds is
closed, we get that it holds for all $\alpha\in\R$, $x\in X_0$,
and $M\in\N$. The claim now follows by letting $M\to \infty$ in
\eqref{E:11},  using our assumption \eqref{E:assumption}, and the
bounded convergence theorem.

{\bf Second step.} Fix $x\in X_0$ and $\e>0$. We claim that for
every $\alpha\in\R$ there exists $N_{\alpha}\in \N$ and open
neighborhood $V_\alpha$ of $\alpha$ ($N_\alpha$ and $V_\alpha$
 also depend on $x$ and $\e$), such that for every
$N>N_{\alpha}$ we have
\begin{equation}\label{E:0}
  \sup_{P\in \R_{k-1}[t], \beta \in V_{\alpha}}\Big|\frac{1}{N}
\sum_{n=1}^N e(n^k\beta+P(n))\cdot f(T^nx)\Big|\leq \e.
\end{equation}
So let $\alpha\in\R$. If $P\in \R_{k-1}[t]$, introducing an error
term $C_{M,N}$ that satisfies $|C_{M,N}|\leq 2 M/N$ we have
\begin{gather}\label{E:3}
 \frac{1}{N} \sum_{n=1}^N  e(n^k\beta+P(n))\cdot f(T^nx)=\\
\notag \frac{1}{[N/M]} \sum_{n=1}^{[N/M]} \frac{1}{M}
\sum_{m=1}^{M}
 e((M n+m)^k\beta+P(Mn+m))\cdot
f(T^{M n+m}x)+C_{M,N}=
\\
\notag A_{M,N}+B_{M,N}+C_{M,N},
\end{gather}
where
$$
A_{M,N}= \frac{1}{[N/M]} \sum_{n=1}^{[N/M]} \frac{1}{M}
\sum_{m=1}^{M} \Big( e\big(m^k\beta+P_{M,n,\beta}(m)\big)
-e\big(m^k\alpha+P_{M,n,\beta}(m)\big)\Big) \cdot f(T^{Mn+m}x)
$$
for some  $P_{M,n,\beta}\in\R_{k-1}[t]$, and
$$
B_{M,N}= \frac{1}{[N/M]} \sum_{n=1}^{[N/M]}\frac{1}{M}
\sum_{m=1}^M e(m^k\alpha+P_{M,n,\beta}(m))\cdot f(T^{Mn+m}x).
$$
 By \eqref{E:key},   there exists
 $M_\alpha\in \N$ such that if  $N$ is large enough we have
$$
|B_{M_\alpha,N}|\leq
\frac{1}{[N/M_\alpha]}\sum_{n=1}^{[N/M_\alpha]} \sup_{P\in
\R_{k-1}[t]}\Big|\frac{1}{M_\alpha} \sum_{m=1}^{M_\alpha}
e(m^k\alpha+P(m))\cdot f(T^{M_\alpha n+m}x)\Big|\leq \e/3.
$$
Choose a neighborhood $V_\alpha$ of $\alpha$ such that
$$
\sup_{\beta\in V_\alpha,1\leq m\leq
M_\alpha}|e(m^k\alpha)-e(m^k\beta)|\leq\e/3.
$$
 Then for  $\beta\in V_\alpha$ we have
$$
|A_{M_\alpha,N}|\leq \frac{1}{[N/M_\alpha]}
\sum_{n=1}^{[N/M_\alpha]}
\frac{1}{M_\alpha}\sum_{m=1}^{M_\alpha}|e(m^k\beta)-e(m^k\alpha)|\leq
\e/3.
$$
Putting all these estimates together we see that  if $N$ is large
 enough, for all $P\in \R_{k-1}[t]$ and $\beta \in
V_\alpha$ the average in \eqref{E:3} is bounded in absolute value
by $2\e/3+M_\alpha/N\leq\e$, proving our claim.

{\bf Third step.} We finish the proof. Notice first that since the
function $e(t)$ is $1$-periodic it suffices to verify \eqref{E:-1}
when the sup is taken over all polynomials with leading term
belonging to the interval $[0,1]$.  Let $\e>0$ and $x\in X_0$. Let
$V_\alpha$ be the open neighborhood of $\alpha\in [0,1]$ and $N_a$
be the positive integer for which \eqref{E:0} holds. Since
$\{V_\alpha\}_{\alpha\in [0,1]}$ is an open cover of the compact
set $[0,1]$ there exists a finite subcover, say
$\{V_{\alpha_i}\}_{i=1,\ldots l}$. Then for every
$N>N_0=\max\{N_{a_1},\ldots, N_{a_l}\}$ we have
$$
\sup_{P\in \R_{k-1}[t],\alpha\in [0,1]}\Big|\frac{1}{N}
\sum_{n=1}^{N} e(n^k\alpha+P(n))\cdot f(T^nx)\Big|\leq \e.
$$
So \eqref{E:-1} holds for every $x\in X_0$,  completing the proof
of the lemma.
\end{proof}
Notice that for totally ergodic systems, the case $k=1$ of
Proposition~\ref{P:key} gives yet another proof of the uniform
Wiener-Wintner theorem (see $(ii)$ in page $1$).

We now combine the argument in \cite{L} with the previous
Proposition to prove Theorem~\ref{T:main}.

\begin{proof}[Proof of Theorem~\ref{T:main}]
The implication $(iii) \Rightarrow (ii)$ is obvious. The
implication $(ii) \Rightarrow (i)$ is easy and was proved in
\cite{L}, but we reprove it for completeness. Let  $g\in
\mathcal{E}_k(T)$, and $f\in L^1(\mu)$ that satisfies condition
$(ii)$ of Theorem~\ref{T:main}. It is easy to check that
\begin{equation}\label{e:0}
g(T^{n}x)=e(P_{x}(n))
\end{equation}
for some $P_{x}\in\R_k[t]$, and  by assumption for a.e. $x\in X$
we have for every $P\in \R_k[t]$ that
\begin{equation}\label{e:1}
\lim_{N\to\infty} \frac{1}{N} \sum_{n=1}^N  e(P(n))\cdot
f(T^nx)=0.
\end{equation}
 By the ergodic theorem and \eqref{e:0} we have for a.e. $x\in X$
 that
$$
\int g f\, d\mu=\lim_{N\to\infty}\frac{1}{N}\sum_{n=1}^N
g(T^nx)\cdot f(T^nx)=\lim_{N\to\infty}\frac{1}{N}\sum_{n=1}^N
e(P_{x}(n))\cdot f(T^{n}x),
$$
and this is $0$ by \eqref{e:1}. Hence,  $f$ is orthogonal to $g$,
showing that $f\in\mathcal{E}(T)^\bot$.

So it remains to show the implication $(i)\Rightarrow (iii)$. By
Proposition~\ref{P:reduction} we can restrict ourselves to the
case where $f\in L^\infty(\mu)$ and $\phi(t)=e(t)$.  We proceed in
two steps:

 {\bf First Step.}  Suppose that
$\alpha\in E_0'$, where $E_0'$  is the set of $\alpha\in\R$ such
that $e(m\alpha)\notin \mathcal{E}_0(T)$ for every nonzero integer
$m$. We claim that for every $k\in \N$ and  $f\in L^\infty(\mu)$,
we have for a.e. $x\in X$ (the set of full measure may depend on
$\alpha$) that
\begin{equation}\label{E22} \lim_{N\to\infty} \sup_{P\in
\R_{k-1}[t]}\Big|\frac{1}{N} \sum_{n=1}^N e(n^k\alpha+P(n))\cdot
f(T^nx)\Big|=0.
\end{equation}
We will use induction on $k$. For $k=1$ the statement follows by
applying the ergodic theorem for the system $(X\times
\mathbb{T},\mathcal{X}\times \mathcal{B},\mu\times m,T\times
R_\alpha)$, and the function $g(x,t)=f(x)\cdot e(t)$, where
$\mathcal{B}$ is the Borel $\sigma$-algebra, $m$ is the Haar
measure on $\mathbb{T}$,  and $R_\alpha(t)=t+\alpha$ (the system
is ergodic because $\alpha\in E_0'$). Suppose that the statement
holds for $k-1$, we will show that it holds for $k$.
 Using van der
Corput's Lemma\footnote{This says that if $a(n)$ is a sequence of
complex numbers then for all integers $1\leq H\leq N$ we have
$$
\Big| \frac{1}{N}\sum_{n=1}^Na(n)\Big|^2\leq
\frac{N+H}{N(H+1)}\cdot\frac{1}{N}\sum_{n=1}^N|a(n)|^2+
2\frac{N+H}{N(H+1)^2}\sum_{h=1}^H(H+1-h)\cdot \text{Re}\Big(
\frac{1}{N}\sum_{n=1}^{N-h} a(n+h)\cdot \overline{a(n)}\Big).
$$}
(see \cite{K}) for the sequence $a_x(n)=e(n^k\alpha+P(n))\cdot
f(T^nx)$, we find that for every $x\in X$ and integers $1\leq
H\leq N$ we have
\begin{gather}\label{E:44}
\sup_{P\in \R_{k-1}[t]}\Big|\frac{1}{N} \sum_{n=1}^N
e(n^k\alpha+P(n))\cdot f(T^nx)\Big|^2 \leq \frac{N+H}{N(H+1)}\cdot
\frac{1}{N}\sum_{n=1}^N|f(T^nx)|^2+\\\notag
2\frac{N+H}{N(H+1)^2}\sum_{h=1}^H(H+1-h)\cdot \sup_{Q\in
\R_{k-2}[t]}\Big| \frac{1}{N}\sum_{n=1}^{N-h}
e(khn^{k-1}\alpha+Q(n))\cdot f(T^{n+h}x)\cdot
\overline{f(T^nx)}\Big|.
\end{gather}
We let $N\to\infty$ in \eqref{E:44}. Using the ergodic theorem for
the function $|f|^2$, and the induction hypothesis for the values
$kh\alpha\in E_0'$ and the functions $T^hf\cdot f$, for $h\in\N$,
we find that for every $H\in\N$ we have for a.e. $x\in X$ that
$$
\lim_{N\to\infty}\sup_{P\in \R_{k-1}[t]}\Big|\frac{1}{N}
\sum_{n=1}^N e(n^k\alpha+P(n))\cdot f(T^nx)\Big|^2\leq
\frac{1}{H+1}\int |f(x)|^2\ d\mu.
$$
Letting $H\to\infty$ we get $\eqref{E22}$, proving the claim.

{\bf Second Step.} We now  prove \eqref{E:main}. We will use
induction on $k$.  For $k=1$ the statement is known to be true for
all ergodic systems (see $(i)$ and $(ii)$ in page $1$). Suppose
that the statement holds for $k-1\geq 1$, we will show that it
holds for $k\geq 2$. By Proposition~\ref{P:key} it suffices to
show that for every $\alpha\in\R$, we have for a.e. $x\in X$ that
\begin{equation}\label{E:nonuniform}
\lim_{N\to\infty} \sup_{P\in \R_{k-1}[t]}\Big|\frac{1}{N}
\sum_{n=1}^N e(n^k\alpha+P(n))\cdot f(T^nx)\Big|=0.
\end{equation}
We consider two cases:

If $\alpha\in E_0'$ then we are covered by the first step.

If $\alpha\notin E_0'$ then  $e(m\alpha)\in \mathcal{E}_0(T)$ for
some nonzero integer $m$. Let $e(\gamma(x))$ be an
$e(m\alpha)$-eigenfunction. Consider the system $(X\times
\mathbb{T}^{k},\mathcal{X}\times \mathcal{B}_k,\mu\times
m_k,T_k)$, where $\mathcal{B}_k$ is the Borel $\sigma$-algebra,
$m_k$ is the Haar measure on $\mathbb{T}^k$,  and
$$
T_k(x,t_1,t_2,\ldots,t_k)=(Tx,t_1+\gamma(x)+b,t_2+t_1,\ldots,t_{k-1}+t_{k-2}).
$$
 As was shown in \cite{L}, it is possible to choose
$b\in\mathbb{T}$ such that the resulting system is totally
ergodic. It was also shown there that for $j\in\mathbb{N}$ every
function in $\mathcal{E}_j(T_k)$ is a product of a function in
$\mathcal{E}_{j+1}(T)$ and a character of $\mathbb{T}^k$. It
follows that if $g\colon X\times\T^k\to \mathbb{C}$ is defined by
$$
g(x,t_1,\ldots,t_k)=f(x)\cdot e(k!\, t_k),
$$
then  $g\in\mathcal{E}_{k-1}(T_k)^\bot$. So we can apply the
inductive hypothesis for the system $(X\times
\mathbb{T}^{k},\mathcal{X}\times \mathcal{B}_k,\mu\times m_k,T_k)$
and the function $g$.  We get that for a.e. $(x,t_1,\ldots,
t_k)\in X\times \T^k$ we have
\begin{equation}\label{E:hypothesis}
\lim_{N\to\infty} \sup_{P\in
\R_{k-1}[t]}\Big|\frac{1}{N} \sum_{n=1}^N e(P(n))\cdot
g(T_k^n(x,t_1,\ldots,t_k))\Big|=0.
\end{equation}
If $C^i_n=\binom{n}{i}$ for $i=1,\ldots,k$, notice that
\begin{align*}
g(T_k^n(x,t_1,\ldots,t_k))=& e\big(k!(t_k+C^1_n\,
t_{k-1}+\ldots+C^{k-2}_n\, t_1+ C^{k-1}_n(\gamma(x)+b)+C^k_n
m\alpha
)\big)\cdot f(T^nx),\\
=& e(n^km\alpha+ Q_{t_1,\ldots,t_k,x,b}(n))\cdot f(T^nx),
\end{align*}
where $Q_{t_1,\ldots,t_k,x,b}\in \R_{k-1}[t]$. Hence,
\eqref{E:hypothesis} gives that for a.e. $x\in X$ we have
\begin{equation}\label{E:nonuniform2}
\lim_{N\to\infty} \sup_{P\in \R_{k-1}[t]}\Big|\frac{1}{N}
\sum_{n=1}^N e(n^km\alpha+P(n))\cdot f(T^nx)\Big|=0.
\end{equation}
So in order to get \eqref{E:nonuniform}, it remains  to replace
$m\alpha$ by $\alpha$ in \eqref{E:nonuniform2}. We do this as
follows: Since $e(m\alpha)\in \mathcal{E}_0(T)$ we have
$e(m^2\alpha)\in \mathcal{E}_0(T^m)$, and since $k\geq 2$ this
implies that  $e(m^k\alpha)\in \mathcal{E}_0(T^m)$. An easy
inductive argument (see \cite{L}, page 784) shows that for totally
ergodic systems $\mathcal{E}_k(T^j)=\mathcal{E}_k(T)$ for
$j\in\N$, so $T^rf\in\mathcal{E}_k(T^m)^\bot$ for every $r\in \N$.
Hence,  we can repeat  the previous argument for the ergodic
system $(X,\mathcal{B},\mu,T^m)$, the eigenvalue $e(m^k\alpha)$ of
$T^m$, and the functions $T^rf$, $r=0,\ldots,m-1$. We get that
for a.e. $x\in X$
$$
\lim_{N\to\infty} \sup_{P\in \R_{k-1}[t]}\Big|\frac{1}{N}
\sum_{n=1}^N e(n^km^k\alpha+P(n))\cdot f(T^{nm+r}x)\Big|=0,
$$
for $r=0,1,\ldots, m-1$. This implies that
$$
\lim_{N\to\infty} \sup_{P\in \R_{k-1}[t]}\Big|\frac{1}{N}
\sum_{n=1}^N e\big((nm+r)^k\alpha+P(nm+r)\big)\cdot
f(T^{nm+r}x)\Big|=0,
$$
for $r=0,1,\ldots, m-1$, which   in turn implies that for a.e.
$x\in X$ we have
$$
\lim_{N\to\infty} \sup_{P\in \R_{k-1}[t]}\Big|\frac{1}{N}
\sum_{n=1}^N e(n^k\alpha+P(n))\cdot f(T^{n}x)\Big|=0.
$$
 This proves \eqref{E:nonuniform} and completes the proof of the induction step.
\end{proof}

\section{Two counterexamples.} We  construct an ergodic  system $(X,\mathcal{B},\mu,T)$ for which
$\mathcal{E}_k(T)^\bot\neq \mathcal{E}_k(T^2)^\bot$ for $k\geq 2$,
and for which part (ii) of Theorem~\ref{T:Lesigne} fails.
\begin{proof}[Proof of Theorem~\ref{T:counterexamples}]
(i) It will be clear that for the system we will construct we have
$\mathcal{E}_k(T)=\mathcal{E}_2(T)$ and
$\mathcal{E}_k(T^2)=\mathcal{E}_2(T^2)$ for $k\geq 2$, so  we can
assume that $k=2$. On $X=\Z_2\times\T^2$ with the Haar measure
$\mu$, consider the measure preserving transformation $T\colon  X
\to X$ defined by
$$
T(0,t_1,t_2)=(1,t_1, t_2), \quad
T(1,t_1,t_2)=(0,t_1+\alpha,t_2+t_1)
$$
for some  $\alpha\in\R$  irrational. We can see that $T$ is
ergodic as follows: We have
$$
T^2(i,t_1,t_2)=(i,t_1+\alpha,t_2+t_1),
$$
and since the skew product transformation
$S(t_1,t_2)=(t_1+\alpha,t_2+t_1)$ defined on $\T^2$ with the Haar
measure is known to be ergodic for $\alpha$ irrational, we have
that $T^2$ has two ergodic components, the sets
$X_i=\{i\}\times\T^2$, $i=1,2$. Since neither of these sets is
$T$-invariant, the transformation $T$ is ergodic.

 It is clear that $e(t_2)\in \mathcal{E}_2(T^2)$. We
will show that $e(t_2)\in \mathcal{E}_2(T)^\bot$. First we compute
$\mathcal{E}_1(T)$. Let  $h\in \mathcal{E}_1(T)$, then
\begin{equation}\label{E:sd}
h(T(i,t_1,t_2))=c\, h(i,t_1,t_2)
\end{equation}
for some nonzero $c\in\C$.  A standard Fourier series argument
shows that $h$ does not depend on $t_2$, so
$h(i,t_1,t_2)=h_1(i,t_1)$. Then \eqref{E:sd} takes the form
$$
h_1(1,t_1)=c\, h_1(0,t_1), \quad  h_1(0,t_1+\alpha)=c\,
h_1(1,t_1).
$$
Substituting the first equation into the second gives
$$
h_1(0,t_1+\alpha)=c^{2} h_1(0,t_1).
$$
 This easily implies that for some $m\in\Z$ and nonzero $c_1\in\C$
we have
 \begin{equation}\label{E:E_1}
h(0,t_1,t_2)=c_1\, e(mt_1), \quad h(1,t_1,t_2)=c_2\, e(mt_1),
\end{equation}
where $c_2=\pm c_1\, e(m\alpha/2)$.

Next we show that $e(t_2)\in \mathcal{E}_2(T)^\bot$. Let
$f\in\mathcal{E}_2(T)$,  then there exists $h\in\mathcal{E}_1(T)$
such that
\begin{equation}\label{E:E_2}
f(T(i,t_1,t_2))=h(i,t_1,t_2)\cdot f(i,t_1,t_2).
\end{equation}
A standard  Fourier series argument  gives that
\begin{equation}\label{E:form}
f(i,t_1,t_2)=e(lt_2)\cdot g(i,t_1),
\end{equation}
for some $l\in\Z$. Substituting this into \eqref{E:E_2} and using
\eqref{E:E_1} gives
\begin{gather*}
 g(1,t_1)= c_1\, e(mt_1)\cdot g(0,t_1),\\
  g(0,t_1+\alpha)\cdot e(lt_2)=c_2\, t e(mt_1)\cdot g(1,t_1),
\end{gather*}
for some $m\in\Z$ and nonzero $c_1,c_2\in\C$. Substituting the
first equation into the second gives
$$
g(0,t_1+\alpha)= c_3\,e((2m-l)t_1)\cdot g(0,t_1)
$$
where $c_3=c_1c_2\neq 0$. The last equation has a solution only if
$l=2m$. Combining this with \eqref{E:form} gives that every
$f\in\mathcal{E}_2(T)$  has the form
$$
f(i,t_1,t_2)=e(2mt_1)\cdot g(i,t_1)
$$
for some $m\in\Z$. This shows that $e(t_2)\in
\mathcal{E}_2(T)^\bot$ and completes the proof of the first part
of Theorem~\ref{T:counterexamples}.

(ii) It suffices to construct an ergodic  system
$(X,\mathcal{B},\mu,T)$,   a function $f\in L^\infty(\mu)$ with
$f\in \mathcal{E}_2(T)^\bot$,
 and a set $X_1\subset X$ with
$\mu(X_1)>0$, and such that for every $x\in X_1$ there exists
$P_x\in \R_2[t]$ satisfying
\begin{equation}\label{E:counter}
\lim_{N\to\infty} \frac{1}{N} \sum_{n=1}^N e(P_x(n))\cdot f(T^nx)
\neq 0.
\end{equation}

 We use the measure preserving system constructed in (i). We
showed before that $f(i,t_1,t_2)=e(t_2)\in\mathcal{E}_2(T)^\bot$.
Notice that for $n\in\N$ we have
\begin{gather*}
T^{2n}(1,t_1,t_2)=(1,t_1+n\alpha,
t_2+nt_1+C_n^2\alpha),\\
T^{2n+1}(1,t_1,t_2)=(0,t_1+(n+1)\alpha,
t_2+(n+1)t_1+C_{n+1}^2\alpha),
\end{gather*}
where $C_n^2=\binom{n}{2}$.
 Let $X_1= \{1\}\times \T^2$, then $\mu(X_1)=1/2$. If
$x=(1,t_1,t_2)\in X_1$ and $P_x\in\R_2[t]$ is such that
$$
P_x(2t)=e\big(-t_2-tt_1-\frac{t(t-1)}{2}\alpha\big),
$$
then for every $n\in\N$
$$
e(P_x(2n))\cdot f(T^{2n}(1,t_1,t_2))=1, \quad e(P_x(2n+1))\cdot
f(T^{2n+1}(1,t_1,t_2))=c_{t_1}e(n\alpha/2)
$$
where $c_{t_1}=e( t_1/2+\alpha/8)$. It follows that
\eqref{E:counter} holds for $x\in X_1$, completing the proof of
the second part of Theorem~\ref{T:counterexamples}.
\end{proof}

\end{document}